\documentclass[twoside]{ima} 
\usepackage{amsmath}
\usepackage{amssymb}
\usepackage{graphicx}
\usepackage{algorithm}
\usepackage{algorithmic}

\newcommand{\cone}{\mathbf{cone}}
\newcommand{\R}{{\mathbb R}}

\newcommand{\N}{\mathbb{N}}

\newcommand{\F}{\mathbb{F}}

\newcommand{\FC}{{\overline{\mathbb{F}}}}
\newcommand{\KC}{{\overline{\mathbb{K}}}}
\newcommand{\K}{\mathbb{K}}

\newcommand{\x}{{\bf x}}
\newcommand{\RD}{{\bar{R}^*}}

\newcommand{\Var}[2]{{\mathcal{V}_{#1}(#2)}}

\newcommand{\genR}[1]{\langle #1 \rangle_R}
\newcommand{\genK}[1]{\langle #1 \rangle_\K}

\begin{document}

\title{Computation with Polynomial Equations and Inequalities arising in Combinatorial Optimization}

\author{Jesus A. De Loera\thanks{Department of Mathematics, University of California at Davis, Davis, CA 95616 (\texttt{deloera@math.ucdavis.edu}); partially supported by NSF and an IBM OCR award}
\and 
Peter N. Malkin\thanks{Department of Mathematics, University of California at Davis, Davis, CA 95616 (\texttt{malkin@math.ucdavis.edu}); partially supported by an IBM OCR award.}
\and Pablo A. Parrilo\thanks{Laboratory for Information and Decision Systems, Department of Electrical Engineering and Computer Science, Massachusetts Institute of Technology, Cambridge, MA 02139 (\texttt{parrilo@mit.edu}); partially supported by AFOSR MURI 2003-07688-1 and NSF FRG DMS-0757207.}
}
\markboth{J.A. DE LOERA, P.N. MALKIN AND P.A. PARRILO}{POLYNOMIALS IN COMBINATORIAL OPTIMIZATION}       

\date{\today}
\maketitle

\begin{abstract}
The purpose of this note is to survey a methodology to solve
systems of polynomial equations and inequalities.  The techniques we
discuss use the algebra of multivariate polynomials with coefficients
over a field to create large-scale linear algebra or semidefinite
programming relaxations of many kinds of feasibility or optimization
questions. We are particularly interested in problems arising in
combinatorial optimization.
\end{abstract}


{\AMSMOS 
90C27, 90C22, 68W05
\endAMSMOS}

\section{Introduction}
\label{sec:intro}

A wide variety of problems in optimization can be easily modeled using
\emph{systems of polynomial equations and inequalities}. Feasibility
and optimization problems translate, either directly or via branching,
into the problem of finding a solution of a system of equations and
inequalities. In this survey paper, we explain how to manipulate such
systems for finding solutions or proving that they do not exist.
Although these techniques work in general, we are particularly
motivated by problems of combinatorial origin. For example, in the
case of graphs, here is how one can think about stable sets,
$k$-colorability and max-cut problems in terms of polynomial
(non-linear) constraints:

\begin{proposition} 
\label{prop:encodings}
Let $G=(V,E)$ be a graph.
\begin{itemize}
\item 
For a given positive integer $k$, consider the following polynomial system:
$$x_i^2 - x_i=0 \; \forall i \in V, \quad x_i x_j=0 
\;\forall (i,j) \in E \; \text{ and } \; \sum_{i\in V}x_i=k.$$
This system is feasible if and only if 
$G$ has a stable set of size $k$.

\item For a positive integer $k$, consider the following polynomial
system of $|V|+|E|$ polynomials equations:
$$x_i^k-1=0 \; \forall i\in V \; \text{ and } \;
\sum_{s=0}^{k-1}x_i^{k-1-s}x_j^s=0\;\forall (i,j)\in E.$$ The graph $G$ is
$k$-colorable if and only if this
system has a complex solution.  Furthermore, when $k$ is odd, $G$ is
$k$-colorable if and only if this system has a common root over
$\FC_2$, the algebraic closure of the finite field with two elements.

\item We can represent the set of cuts of $G$ (i.e., bipartitions on $V$) as the 0-1 incidence vectors
$$SG := \{ \chi^F \,:\, F \subseteq E\,\,\textup{is contained in a cut
  of} \,\, G \} \subseteq \{0,1\}^E.$$ Thus, the max cut problem with
non-negative weights $w_e$ on the edges $e \in E$ is 
$$\textup{max} \{ \sum_{e \in E} w_e x_e \,:\, x \in SG \}.$$
The vectors $\chi^F$ are the solutions of the polynomial system
$$x_e^2 - x_e=0 \,\,\forall \,\,e \in E, \; \text{and} \; \prod_{i \in T} x_i=0
\; \; \forall \,\, T \,\, \textup{an odd cycle in} \,\,G . $$
\end{itemize}
\end{proposition}
There are many other combinatorial problems that can be modeled
concisely by polynomial systems (see \cite{DeLoeraLeeMarguliesOnn08}
and the many references therein). In fact, a given problem can often
be modeled non-linearly in many different ways, and in practice
choosing a ``good'' formulation is critical for an efficient solution.

Given a polynomial system encoding a combinatorial question, we
explain how to use two famous algebraic identities to derive solution
methods. In what follows, let $\K$ denote a field and let $\KC$ denote
the algebraic closure of $\K$.  Let $R=\K[x_1,\ldots,x_n]=\K[x]$ denote
the ring of polynomials in $n$ variables with coefficients over
$\K$. The situation is slightly different depending on whether only
equations are being considered, or if there also inequalities (more
precisely, on whether the underlying field $\K[x]$ is algebraically
closed or formally real):

\begin{enumerate}
\item 
First, suppose that the system contains only the polynomial equations
$f_1(x)=0,f_2(x)=0,\dots,f_s(x)=0$ where $f_1,...,f_s\in\K[x]$.
We explain how to generate a finite sequence of \emph{linear algebra}
systems which terminate with either a solution over $\KC$ of the problem or
provide a certificate of infeasibility. The calculations reduce to matrix
manipulations, mostly rank computations. The techniques we use are a
specialization of prior techniques from computational algebra (see
\cite{Mourrain99,KehreinKreuzer05,KehreinKreuzerRobbiano05,MourrainTrebuchet08}).
As it turns out this technique is particularly effective when the
number of solutions is finite, when $\K$ is a finite field, or when
the system has nice combinatorial information (see
\cite{DeLoeraLeeMarguliesOnn08}). 

\item Second, several authors (see e.g. \cite{Lasserre2000,
Parrilo2003, Laurent2007} and references therein) considered the
solvability (over the reals) of systems of polynomial equations and
inequalities.  It was shown that in this situation there is a way to
set up the feasibility problem
$$\exists x\in\R^n \ s.t. \
f_1(x)=0,\ldots,f_s(x)=0,g_1(x)\ge0,\ldots,g_k(x)\ge0,$$
where $f_1,\ldots,f_s,g_1,\ldots,g_k\in\R[x]$, as a sequence of
semidefinite programs terminating with a feasible solution (see
\cite{Laurent2007}).  Once more, the combinatorial structure can help
in the understanding of the structure of these relaxations, as is
well-known from the case of stable sets \cite{Lovasz1994} and max-cut
\cite{Laurent2004}.  In recent work, Gouveia et al.\
\cite{GouveiaParriloThomas2008,newgouveiaetal2009} considered a
sequence of semidefinite relaxations of the convex hull of real
solutions to an arbitrary combinatorial polynomial system. They called
these approximations {\em theta bodies} because for stable sets of
graphs the first theta body in this hierarchy is exactly Lov{\'a}sz's
theta body of a graph \cite{Lovasz1994}.
\end{enumerate}

The common central idea to both of the relaxations procedures
described above is to use the right \emph{infeasibility certificates}
or \emph{theorems of alternative}. Just as Farkas' lemma is a
centerpiece for the development of Linear Programming, here the key
point is that the infeasibility of polynomial systems can
\emph{always} be certified by particular algebraic identities (on
non-linear polynomials). To find these infeasibility certificates we
rely either on \emph{linear algebra} or \emph{semidefinite
programming} (for a quick overview of semidefinite programming see
\cite{VandenbergheBoyd1996}).

We now state the necessary notation and algebraic concepts that
justify our approach. For a detailed introduction we recommend the
books
\cite{CoxLittleOShea92,CoxLittleOShea05,BochnakCosteRoy1998,Marshall}.
We denote the monomials in the polynomial ring
$R=\K[x_1,\ldots,x_n]=\K[x]$ as
$x^\alpha:=x_1^{\alpha_1}x_2^{\alpha_2}\cdots x_n^{\alpha_n}$ for
$\alpha\in\N^n$.  The degree of $x^\alpha$ is
$\deg(x^\alpha):=|\alpha|:=\sum_{i=1}^n\alpha_i$.  The degree of a polynomial
$f=\sum_{\alpha\in\N^n} f_{\alpha} x^\alpha$, written $\deg(f)$, is the maximum
degree of $x^\alpha$ where $f_\alpha\ne0$ for $\alpha\in\N^n$.  Given a set of
polynomials $F\subset R$, we write $\deg(F)$ for the maximum degree of
the polynomials in $F$.  The variety of $F$ over $\K$, written
$\Var{\K}{F}$, is the set of common zeros of polynomials in $F$ in
$\K^n$, that is, $\Var{\K}{F}:=\{v\in\K^n: f(v)=0\;\forall f\in I\}$.
Also, $\Var{\KC}{F}$, the variety of $F$ over $\KC$, is the set of
common zeros of $F$ in $\KC^n$. Note that in combinatorial problems,
the \emph{variety} of a polynomial system typically has finitely many
solutions (e.g., colorings, cuts, stable sets, etc.).

Given a set of polynomials $F:=\{f_1, \ldots, f_m\}\subseteq R=\K[x]$,
we define the \emph{ideal} of $F$ as
\[ \genR{F}:=\genR{f_1,\ldots, f_m} := \left\{ \sum_{i=1}^m \beta_i f_i \;|
\;\beta_1,\ldots,\beta_m \in \K[x] \right\}.  \] For an ideal
$I\subseteq R$, when $\Var{\KC}{I}$ is finite, the ideal is called
\emph{zero-dimensional} (this is the case for all of the applications
considered here). An ideal $I\subseteq R$ is \emph{radical} if
$f^k\in I$ for some positive integer $k$ implies $f\in I$.
We denote by $\sqrt{I}$ the ideal of all polynomials $f\in R$
such that $f^k\in I$ for some positive integer $k$.
The ideal $\sqrt{I}$ is necessarily radical and it is called the radical ideal
of $I$. Note that $I$ is radical if and only if $I=\sqrt{I}$.

To study of varieties over a non-algebraically closed field like
${\mathbb R}$ requires extra structure.  Given a set of real
polynomials $G:=\{g_1, \ldots, g_m\}\subseteq \R[x]$, we define the
\emph{cone} of $G$ as
\[
\begin{split}
\cone(G) := 
\{
g \, | \, g = s_0 + \sum_{\{i\}} s_i g_i
+ \sum_{\{i,j\}} s_{ij} g_i g_j + \sum_{\{i,j,k\}} s_{ijk} g_i g_j g_k+ \cdots 
\},
\end{split}
\]
where each term in the sum is a square-free product of the polynomials
$g_i$, with a coefficient $s_\alpha \in \R[x]$ that is a sums of
squares. The sum is finite, with a total of $2^m-1$ terms,
corresponding to the nonempty subsets of $\{g_1,\ldots,g_m\}$.

The notions of \emph{ideal} and \emph{cone} are standard in real
algebraic geometry, but they also have inherent convex geometry:
Ideals are affine sets and cones are closed under convex combinations
and non-negative scalings, i.e., they are actually cones in the convex
geometry sense.  Ideals and cones are used for deriving new
\emph{valid constraints}, which are logical consequences of the given
constraints. For example, notice that by construction, every
polynomial in $\genR{f_1,\dots,f_m}$ vanishes in the solution set of
the system $f_1(x)=0,\dots,f_m(x)=0$ over the algebraic closure of
$\K$.  Similarly, every element of $\cone(g_i)$ is clearly
non-negative on the feasible set of $g_1(x) \geq 0,\dots,g_m(x)\ge0$.

It is well-known that optimization algorithms are intimately tied to
the development of feasibility certificates.  For example, the simplex
method is closely related to Farkas' lemma.  Our starting point is a
generalization of this famous principle.  We start with a description
of two powerful infeasibility certificates for polynomial systems
which generalizes the classical ones for linear optimization.  First,
recall from elementary linear algebra the ``Fredholm alternative theorem''
(e.g., \cite{SchrijverLIP}):

\begin{theorem}[Fredholm's alternative]
\label{theo:range}
Given a matrix $A\in\K^{m\times n}$ and a vector $b\in\K^m$,
\[
\nexists\,x\in\K^n \text{ s.t. } A x + b = 0 \; \Leftrightarrow \;
\exists \, \mu\in\K^m \text{ s.t. }  \mu^T A = 0, \; \mu^T b = 1. 
\]   
\end{theorem}
It turns out that there are much stronger versions for general
polynomials, which unfortunately does not seem to be widely known
among optimizers (for more details see e.g., \cite{CoxLittleOShea92}).

\begin{theorem}[Hilbert's Nullstellensatz]
\label{theo:HN}
Let $F:=\{f_1,\ldots,f_m\}\subseteq\K[x]$. Then,
\[
\nexists\,x\in\KC^n \text{ s.t. } f_1(x)=0,...,f_s(x)=0
\Leftrightarrow 1 \in \genR{F}.
\]   
\end{theorem}

Note that $1\in\genR{F}$ means that there exist polynomials
$\beta_1,\ldots,\beta_m\in\K[x]$ such that $1=\sum_{i=1}^m\beta_if_i$.
Note that Fredholm's alternative theorem is simply a special case of
Hilbert's Nullstellensatz where all the polynomials are linear and
the $\beta_i$'s are constant.

Now, the two theorems above deal only with the case of
\emph{equations}. The inclusion of inequalities in the problem
formulation poses additional algebraic challenges because we need to
take into account special properties of the \emph{reals}.  Consider
first the case of \emph{linear} inequalities where linear programming
duality provides the following characterization:

\begin{theorem}[Farkas' lemma]
\label{theo:Farkas}
Let $A\in\R^{m\times n}$, $b\in\R^m$, $C\in\R^{k\times n}$, and $d\in\R^k$.
\begin{gather*}
\nexists\, x\in\R^n \text{ s.t. } A x + b = 0, C x + d \geq 0 \\
\Updownarrow \\
\exists \, \lambda \in\R^m_+,\, \exists \, \mu \in\R^k\; \text{ s.t. } 
\mu^T A + \lambda^T C = 0, \mu^T b + \lambda^T d = -1.
\end{gather*}
\end{theorem}
Again, although not widely known in optimization, it turns out that similar
certificates do exist for \emph{arbitrary} systems of polynomial
equations and inequalities over the reals. The result essentially
appears in this form in \cite{BochnakCosteRoy1998}, and is due to
Stengle \cite{Stengle72}.

\begin{theorem}[Positivstellensatz]
\label{thm:psatz}
Let $F:=\{f_1,\ldots,f_m\}\subset \R[x]$ and $G:=\{g_1,\ldots,g_k\}\subset \R[x]$.
\begin{gather*}
\nexists x\in\R^n \text{ s.t. } f_1(x)=0,\ldots,f_m(x)=0, g_1(x)\ge 0,\ldots,g_k(x)\ge0\\
\Updownarrow \\
\exists \, f\in\genR{F}, \exists \, g \in \cone(G) \text{ s.t. } 
f(x) + g(x) = -1
\end{gather*}
\end{theorem}
The theorem states that for every infeasible system of polynomial
equations and inequalities, there exists a simple algebraic identity
that directly certifies the non-existence of real solutions. 

Of course, we are very concerned with the effective practical
computation of the infeasibility certificates. For the sake of
computation and complexity, we must worry about the growth of degrees
of the infeasibility certificates. On the negative side, the degrees
of the certificates are expected to be high (in the worst case) simply
because the NP-hardness of the original combinatorial questions; see
e.g.\  \cite{DeLoeraLeeMarguliesOnn08}.  At the same time, tight
exponential upper bounds have been derived (see e.g.
\cite{Kollar1988}, \cite{GrigorievVorobjov2002} and references
therein).  Nevertheless, for many problems of practical interest, it
is often the case that it is possible to prove infeasibility using
low-degree certificates (see \cite{DeLoeraLeeMalkinMargulies08,
DeloeraHillarMalkinWoo09}).  Even more important is the fact that for
fixed degree of the certificates the calculations can be reduced to
either linear algebra or semidefinite programming. We summarize the
strong analogies between the case of linear equations and inequalities
with high-degree polynomial systems in the following table:

{\small
\setlength{\tabcolsep}{0.5\tabcolsep}
\begin{table*}[ht]
\begin{center}
\begin{tabular}{c|c|c}
Degree$\backslash$Field \quad & Arbitrary & Real \\ \hline 
Linear & \emph{Fredholm Alternative} & \emph{Farkas' Lemma} \\
       & Linear Algebra & Linear Programming \\ \hline 
Polynomial & \emph{Nullstellensatz} & \emph{Positivstellensatz} \\
& Bounded degree Linear Algebra & \ Bounded degree SDP \\
\end{tabular}
\end{center}
\caption{Infeasibility certificates and their associated computational techniques.}
\label{tab:certificates}
\end{table*}
}

It is important to remark that just as in the classical case of linear
programming, the problem of computation of certificates has very
natural primal-dual formulations, with the corresponding primal and
dual variables playing distinct, but well-define roles. For example,
in the case of Fredholm's alternative, the primal variables are the
variables $x_1,\dots,x_n$ while there is a dual variable for each
equation. For Nullstellensatz and Positivstellensatz there is a
similar duality, based on linear duality and semidefinite programming
duality, respectively. In what follows, we use the most intuitive or
convenient set-up and we leave the reader the exercise of transferring
the results to the corresponding dual version.

The remainder of the paper is divided in two main sections: Section
\ref{sec:equations} is a study of the Hilbert Nullstellensatz, for
general fields, used in the solution of systems of equations. In
Section \ref{sec:inequalities}, we survey the use of the
Positivstellensatz in the context of solving systems of equations and
inequalities over the reals. Both sections contain combinatorial
applications that show why these techniques can be of interest in this
setting. The focus of the combinatorial results is understanding those
situations when a constant degree certificate is enough to show
infeasibility.  These are situations when hard combinatorial problems
have polynomial time algorithms and as such provide structural
insight. Finally, in Section~\ref{sec:solutions}, we describe a
methodology, common to both approaches, to recover feasible solutions
of the original combinatorial problem from the outcome of these
relaxations.

To conclude the introduction we include some more notation.  Given a
vector space $W$ over a field $\K$, we write $\dim(W)$ for the
dimension of $W$.  Given vector spaces $U\subseteq W$, we write $W/U$
as the vector space quotient. Recall that $\dim(W/U)=\dim(W)-\dim(U)$.
As a slight abuse of notation, if $U\not\subseteq W$, then we write
$W/U$ when we strictly mean $W/(U\cap W)$, in which case,
$\dim(W/U)=\dim(W)-\dim(U\cap W)$.  Given a set $F\subset R$,
$\genK{F}$ denotes the vector space generated by $F$ over the
field $\K$. Please note the distinction between the vector space
$\genK{F}$ and the ideal $\genR{F}$.

\section{Solving combinatorial systems of equations}
\label{sec:equations}

In this section, we wish to solve a given zero-dimensional system of
polynomial equations $f_1(x)=0,f_2(x)=0,\dots,f_m(x)=0$ where
$f_1,\ldots,f_s\in R$. 
We abbreviate this system as $F(x)=0$ where $F:=\{f_1,\ldots,f_m\}\subset R$.
Here, by solving a system, we mean first determining if $F(x)=0$ is feasible
over $\KC$, the algebraic closure of $\K$, and furthermore finding a solution
(or all solutions) of $F(x)=0$ if feasible. We say that a system is
\emph{combinatorial} when it is defined in terms of combinatorial information such as
graph properties and it has finitely many solutions (if any).  The
literature on polynomial solving is very extensive and it continues to
be an area of active research (see \cite{Stetter04,CoxLittleOShea05, DickensteinEmiris05}
for an overview and background).  

Here we choose to focus on techniques that fit well with optimization
methods.  The main idea is that solving a polynomial system of
equations can be reduced to solving a sequence of linear algebra
problems.
The foundations of this technique can be traced back to
(\cite{Mourrain99,KehreinKreuzer05,KehreinKreuzerRobbiano05,MourrainTrebuchet08}).
Variants of this technique have been applied to stable sets
\cite{DeLoeraLeeMarguliesOnn08, SusanPhd}, vertex coloring
\cite{DeLoeraLeeMalkinMargulies08, SusanPhd}, satisfiability
(see e.g., \cite{CleggEdmondsImpagliazzo1996}) and
cryptography (see for example \cite{CourtoisKlimovPatarinShamir2000}).
This technique is also strongly related to Gr\"obner bases techniques
(see e.g., \cite{KehreinKreuzer05,MourrainTrebuchet08,Stetter04}).

The linear algebra systems of equations have primal and dual
representations in the sense of Fredholm's lemma.  Specifically, in
this survey, the primal approach solves a linear system to find
constant multipliers $\mu\in\K^m$ such that $1=\sum_{i=1}^m\mu_if_i$
providing a certificate of (non-linear) infeasibility.  Then, the dual
approach aims to find a vector $\lambda$ with entries in $\KC$ indexed
by monomials such that $\sum_{\alpha}\lambda_{x^\alpha}f_{i,\alpha} =
0$ for all $i=1,\dots,m$ and $\lambda_1=1$ where
$f_i=\sum_{\alpha}f_{i,\alpha}x^\alpha$ for all $i$.  As we see in
Section \ref{sec:LArelax}, the dual approach amounts to constructing
linear relaxations of the set of feasible solutions.  In Sections
\ref{sec:LAcert} and \ref{sec:LArelax}, we present examples of the
primal and dual approaches respectively.

\subsection{Linear algebra certificates}
\label{sec:LAcert}

Consider the following consequence of Hilbert's Nullstellensatz: If
there exists \emph{constants} $\mu\in\K^m$ such that
$\sum_{i=1}^m\mu_i f_i = 1$, then the polynomial system $F(x)=0$ must
be infeasible.  In other words, if the system $F(x)=0$ is infeasible,
then $1 \in \genK{F}$.
The crucial point here is that determining whether there exists a $\mu\in\K^m$ 
such that $\sum_{i=1}^m\mu_i f_i = 1$ is a linear algebra problem over $\K$.
The equation $\sum_{i=1}^m\mu_if_i=1$ is called a \emph{certificate of infeasibility}
of the polynomial system.

\begin{example}
\label{ex:1}
Consider the following infeasible system in $\R[x_1,x_2,x_3]$:
\begin{equation*}
x_1^2 - 1 = 0, \; 2x_1x_2+x_3=0, \; x_1+x_2=0, \; x_1+x_3=0.
\end{equation*}
Let $F=\{f_1,f_2,f_3,f_4\}$ where
$f_1=x_1^2 - 1 = 0$, $f_2=2x_1x_2+x_3=0$, $f_3=x_1+x_2=0,$ and $f_4=x_1+x_3=0.$
So, we abbreviate the above system as $F(x)=0$.
We can prove that the system $F(x)=0$ is infeasible if we can find $\mu\in\R^4$
satisfying the following:
\begin{align*}
\mu_1f_1 +\mu_2f_2+\mu_3f_3+\mu_4f_4& =1 \\
\Leftrightarrow 
\mu_1(x_1^2-1)+\mu_2(2x_1x_2+x_3)+\mu_3(x_1+x_2)+\mu_4(x_1+x_3)&=1 \\
\Leftrightarrow 
\mu_1x_1^2+2\mu_2x_1x_2+(\mu_2+\mu_4)x_3+\mu_3x_2+
(\mu_3+\mu_4)x_1-\mu_1&=1.
\end{align*}
Then, equating coefficients on the left and right hand sides of the equation above
gives the following linear system of equations:
\begin{align*}
-\mu_1  &= 1 \quad(1),& \mu_3+\mu_4 &= 0 \quad(x_1),& \mu_3 &= 0 \quad(x_2), \\
\mu_3+\mu_4 &= 0 \quad(x_3),& 2\mu_2 &= 0 \quad(x_1x_2),& \mu_1 &= 0 \quad(x_1^2).
\end{align*}
We abbreviate this system as $\mu^TF=1$.
Even though $F(x)=0$ is infeasible, the linear system
$\mu^TF=1$ is infeasible, and so, we have not found a certificate
of infeasibility.
\end{example}

More formally, let $f_i=\sum_{\alpha\in\N^n}f_{i,\alpha} x^\alpha$ for
$i=1,...,m$.  Note that only finitely many $f_{i,\alpha}$ are non-zero.
Then, $\sum_{i=1}^m\mu_i f_i = 1$ if and only if
$\sum_{i=1}^m \mu_if_{i,0} = 1$ and $\sum_{i=1}^m \mu_if_{i,\alpha} = 0$ for
all $\alpha\in\N^n$ where $\alpha\ne0$.
Note that there is one linear equation per monomial appearing in $F$.
We abbreviate this linear system as $\mu^TF=1$ where we consider $F$ as a
matrix whose rows are the coefficient vectors of its polynomials and we
consider the constant polynomial $1$ as the vector of its coefficients (i.e., a unit vector). The columns of $F$ are indexed by monomials with non-zero
coefficients.

We remark that in the special case where $F(x)=0$ is a linear system
of equations, then Fredholm's alternative says that $F(x)=0$ is
infeasible if and only if $\mu^TF=1$ is feasible.

In general, even if $F(x)=0$ is infeasible, $\mu^TF=1$ may not be
feasible as in the above example.  In order to prove infeasibility, we
must add polynomials from $\genR{F}$ to $F$ and try again to find a
$\mu$ such that $\mu^TF=1$.  Hilbert's Nullstellensatz guarantees
that, if $F(x)=0$ is infeasible, there exists a finite set of
polynomials from $\genR{F}$ that we can add to $F$ so that the linear
system $\mu^T F=1$ is feasible.

More precisely, it is enough to add polynomials of the form $x^\alpha
f$ for $x^\alpha$ a monomial and some polynomial $f\in F$.  Why is
this?  If $F(x)=0$ is infeasible, then Hilbert's Nullstellensatz says
$\sum_{i=1}^m \beta_i f_i=1$ for some $\beta_1,\ldots,\beta_m\in R$.
Let $d=\max_i\{\deg(\beta_i)\}$.  Then, if we add to $F$ all
polynomials of the form $x^\alpha f$ where $f\in F$ and
$\deg(x^\alpha)\le d$.  Then, the $\K$-linear span of $F$, that is
$\genK{F}$, contains $\beta_i f_i$ for all $i$, and thus,
$1\in\genK{F}$ or equivalently $\mu^T F'=1$ is feasible (as a linear
algebra problem) where $F'$ denotes the larger polynomial system.

\begin{example}
Consider again the polynomial system $F(x)=0$ from Example \ref{ex:1}.
Here, $\mu^TF=1$ is feasible, so we must thus add redundant polynomial
equations to the system $F(x)=0$.
In particular, we add the following redundant polynomial equations:
$x_2f_1(x)=0$, $x_1f_2(x)=0$, $x_1f_3(x)=0$, and $x_1f_4(x)= 0$.
Let $F':=\{f_1,f_2,f_3,f_4,x_2f_1,x_1f_2,x_1f_3,x_1f_4\}$.

Then, the system $\mu^T F'=1$ is now as follows:
\begin{align*}
-\mu_1  &= 1 \quad(1),& \mu_3+\mu_4 &= 0 \quad(x_1),& \mu_3-\mu_5 &= 0 \quad(x_2), \\
\mu_2+\mu_4 &= 0 \quad(x_3),& 2\mu_2+\mu_7 &= 0 \quad(x_1x_2),& \mu_1+\mu_7+\mu_8 &= 0 \quad(x_1^2), \\
\mu_6 + \mu_8&= 0 \quad(x_1x_3),& \mu_5 +2\mu_6 &= 0\quad(x_1^2x_2).&
\end{align*}
This system is feasible proving that $F(x)=0$ is infeasible.
The solution is $\mu=(-1,-\frac{2}{3},-\frac{2}{3},\frac{2}{3},-\frac{2}{3},
-\frac{1}{3},\frac{4}{3},-\frac{1}{3})$, which gives the following certificate of infeasibility:
\begin{gather}
-f_1 -\frac{2}{3}f_2 -\frac{2}{3}f_3 +\frac{2}{3}f_4
-\frac{2}{3}x_2f_1 +\frac{1}{3}x_1f_2 +\frac{4}{3}x_1f_3 -\frac{1}{3}x_1f_4 = 1. \notag
\end{gather}
\end{example}
Next, we present the dual approach to the one in this section.

\subsection{Linear algebra relaxations}
\label{sec:LArelax}
In optimization, it is quite common to ``linearize'' non-linear
polynomial systems of equations by replacing all monomials in the system 
with new variables giving a system of linear constraints.
Specifically, we can construct a linear algebra relaxation of the solutions of $F(x)=0$
by replacing every monomial $x^\alpha$ in a polynomial equation in $F(x)=0$
with a new variable $\lambda_{x^\alpha}$ thereby giving a system of linear
equations in the new $\lambda$ variables, one variable for each monomial
appearing in $F$.
Readers familiar with relaxation procedures such as Sherali-Adams and
Lov\'asz-Schrijver (see \cite{Laurent2003} and references therein) will see a
lot of similarities, but here we deal only with \emph{equality} constraints.

\begin{example}
Consider the following feasible system in $\R[x_1,x_2,x_3]$:
\begin{equation*}
f_1(x)=x_1^2 - 1 = 0, \quad f_2(x)=2x_1x_2+x_3=0, \quad f_3(x)=x_1+x_2=0.
\end{equation*}
This system has two solutions $(x_1,x_2,x_3)=(1,-1,2)$ and
$(x_1,x_2,x_3)=(-1,1,2)$.
Let $F=\{f_1,f_2,f_3\}$. So, we abbreviate the above system as $F(x)=0$.
We can replace the monomials $1,x_1,x_2,x_3,x_1^2,x_1x_2$ with the variables
$\lambda_{1},\lambda_{x_1},\lambda_{x_2},\lambda_{x_3},\lambda_{x_1^2},\lambda_{x_1x_2}$
respectively.  The system $F(x)=0$ thus gives rise to the
following set of \emph{linear} equations:
\begin{equation}
\label{y_eqns}
\lambda_{x_1^2} - \lambda_{1} = 0, \qquad 2\lambda_{x_1x_2}+\lambda_{x_3}= 0, \qquad \lambda_{x_1}+\lambda_{x_2}=0.
\end{equation}
We abbreviate the above system as $F*\lambda=0$.

Solutions of $F(x)=0$ give solutions of $F*\lambda=0$:
If $x$ is a solution of $F(x)=0$ above, then
setting $\lambda_{1}=1,\lambda_{x_1}=x_1,\lambda_{x_2}=x_2,
\lambda_{x_3}=x_3,\lambda_{x_1^2}=x_1^2,\lambda_{x_1x_2}=x_1x_2$
gives a solution of $F*\lambda=0$.
So, taking $x=(1,-1,1)$, we set $\lambda_{1}=1$, $\lambda_{x_1}=1$,
$\lambda_{x_2}=-1$, $\lambda_{x_3}=2$, $\lambda_{x_1^2}=1$, and
$\lambda_{x_1x_2}=-1$.  Then, we have $F*\lambda=0$.
Thus, the solutions of $F*\lambda=0$ gives a vector space effectively
containing all of the solutions of $F(x)=0$.
Hence, $F*\lambda=0$ gives a linear relaxation of $F(x)=0$.

There are solutions of $F*\lambda=0$ that do not correspond to
solutions of $F(x)=0$ because the linear system $F*\lambda=0$
does not take into account the non-linear constraints that $\lambda_{1}=1$,
$\lambda_{x_1^2} = \lambda_{x_1}^2$ and
$\lambda_{x_1x_2}=\lambda_{x_1}\lambda_{x_2}$; 
For example, $\lambda_{1}=1$, $\lambda_{x_1}=2$,
$\lambda_{x_2}=-2$, $\lambda_{x_3}=-2$, $\lambda_{x_1^2}=1$ and
$\lambda_{x_1x_2}=1$ is a solution of $F*\lambda=0$, but
$x_1=\lambda_{x_1}=2$, $x_2=\lambda_{x_2}=-2$, and
$x_3=\lambda_{x_3}=-2$ is not a solution of $F(x)=0$.
\end{example}

We now formalize the above example construction of a linear system.
We can consider the polynomial ring $R=\K[x_1,\ldots,x_n]$ as an infinite
dimensional vector space over $\K$ where the set of all monomials $x^\alpha$
forms a vector space basis of $R$.
In other words, a polynomial $f=\sum_{\alpha\in\N^n} f_{\alpha} x^\alpha$
can be represented as an infinite sequence $(f_{\alpha})_{\alpha\in\N^n}$ where
only finitely many $f_{\alpha}$ are non-zero.
We define $\RD=\KC[[x_1,\ldots,x_n]]=\KC[[x]]$ as the ring of formal power series in the
variables $x_1,\ldots,x_n$ with coefficients in $\KC$. 
So, the power series $\lambda=\sum_{\alpha\in\N^n} \lambda_{\alpha} x^\alpha$ can be
represented as an infinite sequence $(\lambda_{\alpha})_{\alpha\in\N^n}$.
Note that we do not require that only finitely many $\lambda_\alpha$ are non-zero.
We define the bilinear form $*:R\times \RD \rightarrow \K$ as follows:
given $f=\sum_{\alpha\in\N^n}f_{\alpha}x^\alpha \in R$ and
$\lambda=\sum_{\alpha\in\N^n} \lambda_{\alpha} x^\alpha \in \RD$,
we define $f*\lambda = \sum_{\alpha\in\N^n} f_{\alpha}\lambda_{\alpha}$,
which is always finite since only finitely many $f_\alpha$ are non-zero.
Thus, we define a linear relaxation of $x\in\KC^n,F(x)=0$, written as
$\lambda\in\RD, F*\lambda=0$, as the set of linear equations $f*\lambda=0$ for
all $f\in F$.

Note that, for any polynomial $f\in R$ and any point $v\in\K^n$, we
have $f(v) = f * \lambda_v$ where $\lambda_v = (v^\alpha)_{\alpha\in\N^n}$.
Thus, for any $v\in\KC^n$, $F(v)=0$ if and only if $F*\lambda_v = 0$.
So, the system $F*\lambda=0$ can be considered as a linear relaxation of
the system $F(x)=0$.
As mentioned in the above example, there are solutions of $F*\lambda=0$ that
do not correspond to solutions of $F(x)=0$ because the linear system
$F*\lambda=0$ does not take into account the relationships between the $\lambda$
variables. Specifically, if $\lambda$ corresponded to a solution of $F(x)=0$,
then we must have $\lambda_{x^\alpha}=\lambda_{x^\beta}\lambda_{x^\gamma}$ 
for all monomials $x^\alpha,x^\beta,x^\gamma$ where
$x^\alpha=x^\beta x^\gamma$.
If we added these non-linear constraints to the linear constraints $F*\lambda=0$,
then we would essentially have the original polynomial system $F(x)=0$.

The system $F*\lambda=0$ is always feasible, but the constraint
$\lambda_{1}=1$ also holds for any $\lambda$ that corresponds to a
solution $x$ of $F(x)=0$. Thus, if the inhomogeneous linear system
$F*\lambda=0$, $\lambda_1=1$ is infeasible, then so is the system of
polynomials $F(x)=0$.
\begin{remark}
Crucially, this linear system $F*\lambda=0$, $\lambda_1=1$ is dual to
the linear system $\mu^TF=1$ from the previous section by Fredholm's
alternative meaning that $F*\lambda=0$, $\lambda_1=1$ is infeasible
if and only if $\mu^TF=1$ is feasible.
\end{remark}

There is a fundamental observation we wish to make here:
\emph{adding redundant polynomial equations can lead to a tighter relaxation}.
\begin{example} (Cont.)
Add $x_1f_3(x)=x_1^2+x_1x_2=0$ to the system $F(x)=0$ giving the system
$F'(x)=0$ where $F':=\{f_1,f_2,f_3,x_1f_3\}$.
The system $F'(x)=0$ has the same solutions as $F(x)=0$.
The polynomial equation $x_1f_3(x)=0$ gives rise to a new linear equation
$\lambda_{x_1^2}+\lambda_{x_1x_2} = 0$ giving the following linear system
$F'*\lambda=0$:
\begin{equation}
\label{y1_eqns}
\lambda_{x_1^2} - \lambda_{1} = 0, \; 2\lambda_{x_1x_2}+\lambda_{x_3}= 0,
\; \lambda_{x_1}+\lambda_{x_2}=0, \; \lambda_{x_1^2}+\lambda_{x_1x_2} = 0.
\end{equation}
The dimension of the solution space of the original system $F*\lambda=0$ is three 
if we ignore all $\lambda$ variables that do not appear in the linear system,
or in other words, if we project the solution space onto the $\lambda$ variables
appearing in the system.
However, the dimension of the projected solution space of $F'*\lambda=0$ is two; so,
$F'*\lambda=0$ is a tighter relaxation of $F(x)=0$.
\end{example}

We denote the set of solutions of the linear system $F*\lambda=0$ as
$F^\circ := \{\lambda \in \RD: F*\lambda=0\}$, called the annihilator of $F$,
which is a vector subspace of $\RD$.  The fact that adding redundant equations
leads to a tighter linear relaxation is summarized by the following fact:  For
sets $F\subseteq \tilde{F}\subseteq R$, we have $F^\circ\subseteq {\tilde{F}^\circ}$.

Extending this idea, consider the ideal $I=\genR{F}$, which is the set of all
redundant polynomials given as a polynomial combination of polynomials in $F$,
then $I^\circ$ becomes a finite dimensional vector space where $\dim(I^\circ)$
is precisely the number of solutions of $F(x)=0$ over $\KC$, including
multiplicities, assuming that there are finitely many solutions.
Note that by linear algebra, $I^\circ$ is isomorphic to the vector
space quotient $R/I$ (see e.g., \cite{Stetter04}).
Furthermore, if $I$ is radical, then $\dim(I^\circ)=\dim(R/I)$ is precisely the number
of solutions of $F(x)=0$.
So, there is a direct relationship between the 
number of solutions of a polynomial system and the dimension of the 
solution space of its linear relaxation (see e.g., \cite{CoxLittleOShea05}).

\begin{theorem}\label{thm: orth ideal dim}
Let $I\subseteq R$ be a zero-dimensional ideal.  Then, $\dim(I^\circ)$ is
finite and $\dim(I^\circ)$ is the number of solutions of polynomial system
$I(x)=0$ over $\KC$ including multiplicities, so $|\Var{\KC}{I}| \le
\dim(I^\circ)$ with equality when $I$ is radical.
\end{theorem}

So, if we can compute $\dim(I^\circ)$, then we can determine the
feasibility of $I(x)=0$ over $\KC$. Unfortunately, we cannot compute
$\dim(I^\circ)$ directly. Instead, under some conditions 
(see Theorem \ref{thm: main}), we can compute $\dim(I^\circ)$ by computing the
dimension of $F^\circ$ when projected onto the $\lambda_{x^\alpha}$ variables
where $\deg(x^\alpha)\le\deg(F)$.

\subsection{Nullstellensatz Linear Algebra Algorithm (NulLA)}
We now present an algorithm for determining whether a polynomial system of
equations is infeasible using linear relaxations.
Let $F\subseteq\K[x]$ and again let $F(x)=0$ be the polynomial system 
$f(x)=0$ for all $f\in F$.
We wish to determine whether $F(x)=0$ has a solution over $\KC$.

The idea behind NulLA \cite{DeLoeraLeeMalkinMargulies08} is
straightforward: we check whether the linear system $F*\lambda =0$,
$\lambda_1=1$ is infeasible or equivalently whether $\mu^TF = 1$ is feasible
(i.e., $1\in\genK{F}$) using linear algebra over $\K$ and if not then we add
polynomials from $\genR{F}$ to $F$ and try again.
We add polynomials in the following systematic way: for each polynomial $f\in
F$ and for each variable $x_i$, we add $x_i f$ to $F$.  So, the NulLA algorithm
is as follows: if $F*\lambda=0,\lambda_1=1$ is infeasible, then $F(x)=0$ is
infeasible and stop, otherwise for every variable $x_i$ and every $f\in F$ add
$x_if$ to $F$ and repeat.

In the following, we assume without loss of generality that $F$ is closed under
$\K$-linear combinations, that is $F=\genK{F}$, and thus, $F$ is a vector space
over $\K$.  Note that taking the closure of $F$ under $\K$-linear combinations
does not change the set of solutions of $F(x)=0$ and does
not change the set of solutions of $F*\lambda=0$.
In practice, we must choose a vector space basis of $F$ for
computation, but the point we wish to make is that the choice of basis
is irrelevant.  Moreover, we find that it is more natural to work with
vector spaces and that it leads to a more concise exposition.  Recall
from above that $F*\lambda=0,\lambda_1=1$ is infeasible if and only if
$1\in \genK{F}$, which when $F$ is a vector space,
simplifies to $1\in F$ since $\genK{F}=F$.

For a vector space $F\subset R$, we define
$F^+ := F + \sum^n_{i=1} x_i F$ where $x_iF:= \{x_if:f\in F\}$.  Note that
$F^+$ is also a vector subspace of $R$.  Then, $F^+$ is precisely the linear
span of $F$ and $x_iF$ for all $i=1,\ldots,n$.
So, the NulLA algorithm for vector spaces is as follows (see
Algorithm \ref{alg:NulLA}): if $1 \in F$, then $F(x)=0$ is infeasible
and stop, otherwise set $F:=F^+$ and repeat.  There is an upper bound
on the number of times we need to repeat the above step given by the
Nullstellensatz bound of the system $F(x)=0$.  This follows since
after $d$ iterations of NulLA, the set $F$ contains all linear
combinations of polynomials of the form $x^\alpha f$ where the total
degree $|\alpha|\le d$ and where $f$ was one of the initial
polynomials in $F$.

\begin{algorithm}
\caption{NulLA Algorithm \cite{DeLoeraLeeMalkinMargulies08}}
\begin{algorithmic}
\REQUIRE A finite dimensional vector space $F\subseteq R$ and a Nullstellensatz bound $D$.
\ENSURE \textsc{Feasible}, if $F(x)=0$ is feasible over $\KC$, else \textsc{Infeasible}.
\FOR{$d = 0,1,2,\ldots,D $}
  \STATE If $1 \in F$, then \textbf{return} \textsc{Infeasible}.
  \STATE $F:=F^+$.
\ENDFOR
\STATE \textbf{Return} \textsc{Feasible.}\\
\end{algorithmic}
\label{alg:NulLA}
\end{algorithm}

While theoretically the Nullstellensatz bound limits the number of
iterations, this bound is in general too large to be practically useful (see
\cite{DeLoeraLeeMalkinMargulies08}).  Hence, in practice, NulLA is
most useful for proving infeasibility (see Section
\ref{sec:Experimental Results}).

Next, we discuss improving NulLA by adding redundant polynomials to $F$
in such a way so that $\deg(F)$ does not grow unnecessarily.
The improved algorithm is called the Fixed-Point Nullstellensatz Linear Algebra
(FPNulLA) algorithm (see \cite{DeloeraHillarMalkinWoo09}).
The basic idea behind the FPNulLA algorithm is that, if $1\not\in F$, then
instead of replacing $F$ with $F^+$ and thereby increasing $\deg(F)$, we check
to see whether there are any new polynomials in $F^+$ with degree at most
$\deg(F)$ that were not in $F$ and add them to $F$, and then check again
whether $1\not\in F$.
More formally, if $1\not\in F$, then we replace $F$ with $F^+\cap R_d$ 
where $R_d$ is the set of all polynomials with degree at most $d=\deg(F)$.
We keep replacing $F$ with $F^+\cap R_d$ until either $1\in F$ or
we reach a \emph{fixed point}, $F=F^+\cap R_d$.
This process must terminate.

Note that if we find that $1\in F$ at some stage of FPNulLA this implies
that there exists an infeasibility certificate of the form
$1=\sum_{i=1}^s \beta_i f_i$ where $\beta_1,...,\beta_s\in\K[x]$ and
the polynomials $f_1,...,f_s\in\K[x]$ are a vector space basis
of the original set $F$.

Moreover, we can also improve NulLA by proving that the system $F(x)=0$ is
feasible well before reaching the Nullstellensatz bound as follows. 
When $1\not\in F$ and $F=F^+\cap R_d$, then we could set $F:=F^+$ and $d:=d+1$
and repeat the above process.  However, when we reach the fixed point
$F=F^+\cap R_d$, we can use the following theorem to determine if the system is
feasible and if so how many solutions it has.
First, we introduce some notation.
Let $\pi_d:\RD \rightarrow \bar{R}_d$ be the projection of a power series onto 
a polynomial of degree at most $d$ with coefficients in $\KC$.
Below, we abbreviate $\dim(\pi_d(F^\circ))$ as $\dim_d(F^\circ)$
and similarly $\dim(\pi_{d-1}(F^\circ))$ as $\dim_{d-1}(F^\circ)$.
\begin{theorem}
\label{thm: main}
Let $F\subset R$ be a finite dimensional vector space and let $d=\deg(F)$.
If $F=F^+\cap R_d$ and %$\dim(R_d/F)=\dim(R_{d-1}/F)$
$\dim_d(F^\circ)=\dim_{d-1}(F^\circ)$,
then $\dim(I^\circ)=\dim_d(F^\circ)$ where $I=\genR{F}$. 
\end{theorem}

See \cite{Mourrain99,DeloeraHillarMalkinWoo09} for a proof of Theorem \ref{thm: main}.
Recall from Theorem \ref{thm: orth ideal dim}, that there are $\dim(I^\circ)$
solutions of $F(x)=0$ over $\KC$ including multiplicities where $I=\genR{F}$
and exactly $\dim(I^\circ)$ solutions when $I$ is radical.

There are many equivalent forms of the above theorem that appear in
the literature.
(see e.g., \cite{Mourrain99,ReidZhi08,LasserreLaurentRostalski2009}).
Note that the condition that $F=F^+\cap R_d$ is
equivalent to the condition that
$\dim_d(F^\circ)=\dim_d((F^+)^\circ))$.
Also, since $R_d$ is a vector space and $F\subseteq R_d$ is vector subspace, we
can form the vector space quotient $R_d/F$, which is isomorphic to
$\pi_d(F^\circ)$ (see for example \cite{Stetter04}), and thus,
$\dim_d(F^\circ)=\dim(R_d/F)=\dim(R_d)-\dim(F)$ where $\dim(R_d)= \binom{n+d}{d}$.
Similarly, $\dim(R_{d-1}/F)=\dim_{d-1}(F^\circ)$ and $\dim(R_d/F^+)=\dim_d((F^+)^\circ)$.
Thus, in practice, checking the conditions of Theorem \ref{thm: main} means
computing $\dim(F)$, $\dim(F\cap R_{d-1})$ and $\dim(F^+\cap R_d)$.

We can now present the FPNulLA algorithm
\cite{Mourrain99,DeloeraHillarMalkinWoo09}.
See \cite{Mourrain99} for a proof of termination.

\begin{algorithm}
\caption{FPNulLA Algorithm}
\begin{algorithmic}
\REQUIRE A vector space $F\subset R$.
\ENSURE The number of solutions of $F(x)=0$ over $\KC$ up to multiplicities.
\STATE Let $d=\deg(F)$.
\LOOP
\STATE \textbf{if} $1 \in F$ \textbf{then} Return $0$.
\WHILE{$F\ne F^+\cap R_d$}
  \STATE Set $F:=F^+\cap R_d$.
  \STATE \textbf{if} $1 \in F$ \textbf{then} return $0$.
\ENDWHILE
  \STATE \textbf{if} $\dim_d(F^\circ)=\dim_{d-1}(F^\circ)$ \textbf{then} return $\dim_d(F^\circ)$.
  \STATE $F:=F^+$.
  \STATE $d:=d+1$.
\ENDLOOP
\end{algorithmic}
\label{alg:FPNulLA}
\end{algorithm}

\begin{example}
Consider the following feasible system with polynomials in $\K[x]$ with
$\K=\F_2$.
\[ 1+x+x^2=0, \quad 1+y+y^2=0, \quad x^2+xy+y^2=0.\]
This system has two solutions over $\KC=\FC_2$.
Let $F:=\genK{1+x+x^2, 1+y+y^2, x^2+xy+y^2}$.
Then, $1\not\in F$ and $\deg(F)=2$.
Now, 
\begin{align*}
F^+=F&+xF+yF\\
 =F&+\genK{x+x^2+x^3, x+xy+xy^2, x^3+x^2y+xy^2} \\
&+ \genK{y+xy+x^2y, y+y^2+y^3, x^2y+xy^2+y^3}
\end{align*}
Then, $F^+\cap R_2=\genK{1+x+x^2, 1+y+y^2, x^2+xy+y^2, 1+x+y}.$
So, $F\ne F^+\cap R_2$.
Next, let $F:=F^+\cap R_2$. Then, $F=F^+\cap R_2$.
Moreover, $\dim_2(F^\circ)=2$ and $\dim_1(F^\circ)=2$.
Therefore, the system is feasible.
\end{example}

\subsection{Experimental results}
\label{sec:Experimental Results}
In this section, we summarize experimental results for graph 3-coloring
from \cite{DeloeraHillarMalkinWoo09},
which illustrate the practical performance of the NulLA and FPNulLA algorithms.
For further and more detailed results, see
\cite{DeLoeraLeeMalkinMargulies08,SusanPhd,DeloeraHillarMalkinWoo09}.
Experimentally, for graph 3-coloring, NulLA and FPNulLA are well-suited to
proving infeasibility, that is, that no 3-coloring exists.
The polynomials encoding of 3-coloring that is used here
is over $\F_2$ (see Proposition \ref{prop:encodings}) and thus any linear
algebra operations are very fast.
However, even though in theory NulLA and FPNulLA can determine
feasibility, for the experiments described below NulLA and FPNulLA
were not able to prove feasibility in practice.

We refer to the number of iterations that NulLA takes to solve a given
system of equations as the \emph{NulLA degree} of the system.
Similarly to the NulLA degree, we refer to the number of outer iterations that
FPNulLA takes to the system as the \emph{FPNulLA degree} of the system.
We can consider the NulLA degree and the FPNulLA degree as measures of the
\emph{hardness} of proving infeasibility of the system.
In this section, we present experimental evidence that
the NulLA degree of an infeasible combinatorial system is a \emph{good}
measure of the \emph{hardness} of proving infeasibility of the system.
Similarly, we present experimental evidence
(see also \cite{CleggEdmondsImpagliazzo1996} for theoretical evidence)
suggesting that the FPNulLA degree is also a \emph{good} measure of the
\emph{hardness} of a problem and an even \emph{better} measure than the NulLA
degree.

Here, we are interested in the percentage of randomly generated graphs
whose polynomial system encoding has a NulLA degree of one or a
FPNulLA degree of one.  The $G(n,p)$ model \cite{Gilbert1959} is used for
generating random graphs where $n$ is the number of vertices and $p$ is the
probability that an edge is included between any two vertices.  Also, without
loss of generality, for a slightly smaller polynomial encoding, the color of
one of the vertices of each randomly generated graph was fixed.

The experimental results are presented in Figure \ref{fig:fpnulla
comp} (taken from \cite{DeloeraHillarMalkinWoo09}),
which plots the percentage of 1000 random graphs in $G(100,p)$ that
were proven infeasible with a NulLA degree of one, with a FPNulLA degree of
one, or with an exact method versus the $p$ value.
The exact method used was to model graph 3-coloring as a Boolean
satisfiability problem \cite{VanGelder2008} and then use the program
\texttt{zchaff} \cite{zchaff} to solve the satisfiability problem.

\begin{figure}[ht]
\begin{center}
\includegraphics[scale=0.8]{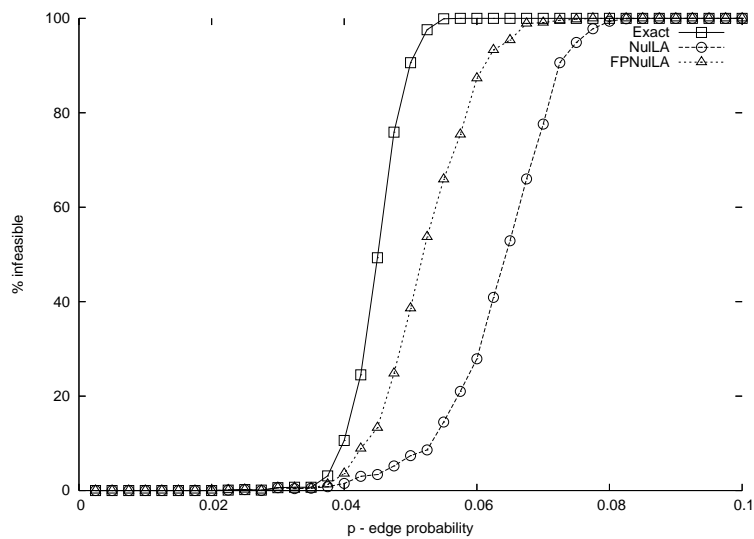}
\caption{Non-3-colorable graphs with NulLA or FPNulLA degree of 1}
\label{fig:fpnulla comp}
\end{center}
\end{figure}

It is well-known that there is a distinct phase transition from
feasibility to infeasibility for graph 3-coloring, and it is at this
phase transition that graphs exists for which it is difficult on
average to prove infeasibility or feasibility (see
\cite{HoggWilliams1994}).  Observe that the infeasibility curve for
NulLA resembles that of the exact infeasibility curve and that the
infeasibility curve for FPNulLA also resembles the infeasibility curve
and clearly dominates the infeasibility curve for NulLA.
These results support that statement that the NulLA degree or FPNulLA degree is
a reasonable measure of the hardness of proving infeasibility since those
graphs that require a higher degree than one are located near the phase
transition.

\subsection{Application: the structure of non-3-colorable graphs}
\label{Application: Non-3-colorable graphs}

For a given class of combinatorial system of equations, it is of
interest to understand the growth of the NulLA degree or FPNulLA
degree.  For some fixed degree, it is also interesting to characterize
which graphs can be proved at that degree to lack a certain property.
In this section, we state a combinatorial characterization of those
graphs whose combinatorial system of equations encoding 3-colorability
has a NulLA degree of one and recall bounds for the NulLA degree (see
\cite{SusanPhd}):

\begin{theorem}
The NulLA degree for a polynomial encoding over $\F_2$ of the
3-colorability of a graph with $n$ vertices with no 3-coloring is at least one
and at most $2n$.  Moreover, in the case of a non-3-colorable graph containing
an odd-wheel or a 4-clique as a subgraph, the NulLA degree is exactly one.
\end{theorem}

\begin{figure}[ht]
\label{fig_3and4cycles}
\begin{center}
\includegraphics[scale=0.3]{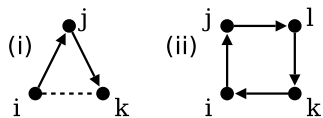}
\end{center}
\caption{(i) oriented partial 3-cycle and (ii) an oriented chordless 4-cycle}
\end{figure}

Now we look at those non-3-colorable graphs that have a degree one NulLA degree.
Let $A$ denote the set of all possible directed edges or arcs in the graph $G$.
We are interested in two types of substructures of the graph $G$:
oriented partial-3-cycles and oriented chordless 4-cycles
(see Figure \ref{fig_3and4cycles}).
An \textbf{oriented partial-3-cycle} is a set of two arcs of a 3-cycle,
that is, a set $\{(i,j),(j,k)\}$ also denoted $(i,j,k)$ where
$(i,j),(j,k),(k,i)\in A$.  
An \textbf{oriented chordless 4-cycle} is a set of four arcs
$\{(i,j),(j,l),(l,k),(k,i)\}$ also denoted $(i,j,k,l)$ where
$(i,j),(j,l),(l,k),(k,i)\in A$ and $(j,k),(i,l)\not\in A$.

Now, we can state a sufficient condition for non-3-colorability
\cite{DeloeraHillarMalkinWoo09}.
This sufficient condition is satisfied if and only if the combinatorial system
encoding 3-coloring has a NulLA degree of one, which is proved in
\cite{DeloeraHillarMalkinWoo09}.
\begin{theorem}
\label{thm:non3color}
The graph $G$ is not 3-colorable if there exists a set $C$ of oriented partial
3-cycles and oriented chordless 4-cycles such that
\begin{enumerate}
\item $|C_{(i,j)}|+|C_{(j,i)}| \equiv 0 \pmod{2}$ for all $(i,j)\in E$ and
\item $\sum_{(i,j)\in A,i<j} |C_{(i,j)}| \equiv 1 \pmod{2}$
\end{enumerate}
where $|C_{(i,j)}|$ denotes the number of cycles in $C$ 
(either 3-cycles or 4-cycles) in which the arc $(i,j)\in A$ appears.
\end{theorem}

Condition 1 in Lemma \ref{thm:non3color} means that every undirected edge of
$G$ is covered by an even number of directed edges from cycles in $C$ (ignoring
orientation).
Condition 2 in Lemma \ref{thm:non3color} means that, given any orientation of
$G$, the total number of times the arcs in that orientation appear in the
cycles of $C$ is odd.
The particular orientation we use in Lemma \ref{thm:non3color} is the
orientation given by the set of arcs $\{(i,j)\in A: i<j\}$, but
the particular orientation we use for Condition 2 is irrelevant 
(see \cite{DeloeraHillarMalkinWoo09}).

\begin{example}
Consider the Gr\"otzsch graph (Mycielski 4) in
Figure~\ref{fig_groetzsch}, which has no 3-coloring.  It contains no
3-cycles.  Now, consider the following set of oriented chordless
4-cycles, which we show gives a certificate of non-3-colorability by
Lemma~\ref{thm:non3color}.
\begin{align*}
C:=\{ &(1,2,3,7), (2,3,4,8), (3,4,5,9),
(4,5,1,10), (1,10,11,7), \\ 
&(2,6,11,8),(3,7,11,9),(4,8,11,10),(5,9,11,6)\}.
\end{align*}
Figure \ref{fig_groetzsch} illustrates the edge directions for
the 4-cycles of $C$.
Each undirected edge of the graph is contained in exactly two
4-cycles, so $C$ satisfies Condition 1 of Lemma \ref{thm:non3color}.
Now, 
$$|C_{(6,11)}|=|C_{(7,11)}|=|C_{(8,11)}|=|C_{(9,11)}|=|C_{(10,11)}|=1,$$
and $|C_{(i,j)}\equiv0\pmod{2}$ for all other arcs $(i,j)\in A$ where $i<j$.
Thus, $$\sum_{(i,j)\in A,i<j} |C_{(i,j)}| \equiv 1 \pmod{2},$$
so Condition 2 is satisfied, and therefore, the graph has no 3-coloring.
\begin{figure}[ht]
\begin{center}
% l b r t
\includegraphics[scale=0.2, trim=0 0 25 0]{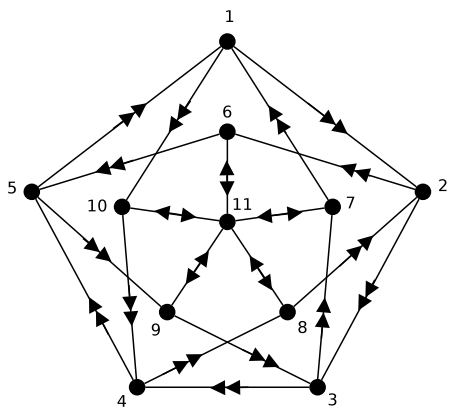}
\caption{Gr\"otzsch graph}
\label{fig_groetzsch}
\end{center}
\end{figure}
\end{example}

\section{Adding polynomial inequalities}
\label{sec:inequalities}

Up until this point we have worked over arbitrary fields (with special
attention to finite fields due to their fast and exact computation),
where the only allowable constraints were equations. Now we turn our
attention to the real case (i.e. $\K = \R$), where we have the
additional possibility of specifying \emph{inequalities} (more
generally, one can work over \emph{ordered} or \emph{formally real}
fields). In this case, following the terminology of real algebraic
geometry, we call the solution set of a system of polynomial equations
and inequalities a {\em basic semialgebraic set}. Note that convex
polyhedra correspond to the particular case where all the constraint
polynomials have degree one. As we have seen earlier in the
Positivstellensatz (Theorem~\ref{thm:psatz} above), the emptiness of a
basic semialgebraic set can be certified through an algebraic identity
involving sum of squares of polynomials.

The connection between sum of squares decompositions of polynomials
and convex optimization can be traced back to the work of N.~Z.~Shor
\cite{Shor}.  His work went relatively unnoticed for several years,
until several authors, including Lasserre, Nesterov, and Parrilo,
observed around the year 2000 that the existence of sum of squares
decompositions and the search for infeasibility certificates for a
semialgebraic set can be addressed via a sequence of semidefinite
programs relaxations \cite{Lasserre2000, Parrilo2000, Parrilo2003,
NesterovSquared}. The first part of this section will be a short
description of the connections between sums of squares and
semidefinite programming, and how the Positivstellensatz allows, in a
analogous way to what was presented in Section~\ref{sec:equations} for
the Nullstellensatz, for a systematic way to formulate these
semidefinite relaxations.

A very central preoccupation of combinatorial optimizers has been the
understanding of the facets that describe the integer hull (normally
binary) of a combinatorial problem.  As we will see in the last part
of this survey, one can recover quite a bit of information about the
integer hull of combinatorial problems from a sequence combinatorially
controlled SDPs. This kind of approach was pioneered in the
lift-and-project method of Balas, Ceria and Cornu\'ejols
\cite{BalasCeriaCornuejols1993}, the matrix-cut method of Lov\'asz and
Schrijver \cite{LovaszSchrijver1991} and the linearization technique
of Sherali-Adams \cite{SheraliAdams1990}. Here we try to present more
recent developments (see \cite{Laurent2009} and references therein for
a very extensive survey).

\subsection{Sums of squares, SDP, and feasibility of semialgebraic sets}

A multivariate polynomial $p(x)$ is a \emph{sum of squares} (SOS for
short) if it can be written as a sum of squares of other polynomials,
i.e.,
\[
p(x) = \sum_i q_i^2(x), \quad q_i(x) \in \R[x].
\]
If $p(x)$ is SOS, then clearly $p(x) \geq 0$ for all $x\in\R^n$. 
\begin{example}
The polynomial $p(x_1,x_2) = x_1^2-x_1 x_2^2+x_2^4+1$ is SOS. Among
infinitely many others, it has the following decompositions:
\begin{eqnarray*}
p(x_1,x_2) &=& \frac{3}{4} (x_1-x_2^2)^2 + \frac{1}{4} (x_1+x_2^2)^2+ 1 \\
&=& 
\frac{1}{9}( 3 -  x_2^2)^2 + 
\frac{2}{3} x_2^2 +
\frac{1}{288} (9 x_1 - 16 x_2^2 )^2 +
\frac{23}{32} x_1^2.
\end{eqnarray*}
\label{ex:poly}
\end{example}
The sum of squares condition is a quite natural sufficient test for
polynomial non-negativity. Thus instead of asking whether even degree
polynomials are non-negative we ask the easier question whether they
are sums of squares. More importantly, as we shall see, the existence
of a sum of squares decomposition can be decided via semidefinite
programming.

\begin{theorem} \label{sostosdp}
A polynomial $p(x)$ is SOS if and only if $p(x) = z^T Q z$, where $z$
is a vector of monomials in the $x_i$ variables, and $Q$ is a
symmetric positive semidefinite matrix.
\end{theorem}

By the theorem above, every SOS polynomial can be written as a
quadratic form in a set of monomials, with the corresponding matrix
being positive semidefinite. The vector of monomials $z$ in general
depends on the degree and sparsity pattern of $p(x)$. If $p(x)$ has
$n$ variables and total degree $2d$, then $z$ can always be chosen as
a subset of the set of monomials of degree less than or equal to $d$,
which has cardinality $\binom{n+d}{d}$.

\begin{example}
Consider again the polynomial from Example~\ref{ex:poly}. It has the
representation
\[
p(x_1,x_2) =\frac{1}{6} 
\left [\begin {array}{c} 1\\x_2\\
x_2^{2}\\x_1\end {array}\right ]^T
\left [\begin {array}{rrrr} 
6&0&-2&0\\
0&4&0&0\\
-2&0&6&-3\\
0&0&-3&6
\end {array}\right ]
\left [\begin {array}{c} 1\\x_2\\
x_2^{2}\\x_1\end {array}\right ],
\]
and the matrix in the expression above is positive semidefinite.
\end{example}

In the representation $f(x) = z^T Q z$, for the right- and left-hand
sides to be identical, all the coefficients of the corresponding
polynomials should be equal. Since $Q$ is simultaneously constrained
by linear equations and a positive semidefiniteness condition, the
problem can be easily seen to be directly equivalent to an
semidefinite programming feasibility problem in the standard primal
form.

Now we describe an algorithm, and illustrate it with an example, on
how we can use SDPs to decide the feasibility of a system of
polynomial inequalities. Exactly as we did for the Nullstellensatz
case, we can look for the existence of a Positivstellensatz
certificate of bounded degree $D$. Once we assume that the degree $D$
is fixed we can apply Theorem~\ref{sostosdp} and obtain a
reformulation as a semidefinite programming problem. We formalize this
description in the following algorithm:

\begin{algorithm}
\caption{Bounded degree Positivstellensatz \cite{Parrilo2000,Parrilo2003}}
\begin{algorithmic}
\REQUIRE A polynomial system $\{f_i(x) =0, g_i(x) \geq 0\}$ and a Positivstellensatz bound $D$.
\ENSURE \textsc{Feasible}, if $\{f_i(x)=0, g_i(x) \geq 0\}$ is feasible over $\R$, 
else \textsc{Infeasible}.
\FOR{$d = 0,1,2,\ldots,D $}
  \STATE If there exist $\beta_i, \, s_\alpha \in \R[x]$ such that $-1=\sum_i \beta_i f_i + \sum_{\alpha \in \{0,1\}^n} s_\alpha g^\alpha$, with $s_\alpha$ SOS, 
$\text{deg}(\beta_i f_i) \leq d$, $\text{deg}(s_\alpha g^\alpha) \leq d$ then \textbf{return} \textsc{Infeasible}.
  \STATE $d:=d+1$.
\ENDFOR
\STATE \textbf{Return} \textsc{Feasible.}\\
\end{algorithmic}
\label{alg:psatz}
\end{algorithm}
Notice that the membership test in the main loop of the algorithm is,
by the results described at the beginning of this section, equivalent
to a finite-sized semidefinite program. Similarly to the
Nullstellensatz case, the number of iterations (i.e., the degree of
the certificates) serves as a quantitative measure of the hardness in
proving infeasibility of the system. As we will describe in more
detail in Section~\ref{sec:chull}, in several situations one can give
further refined characterization on these degrees.

\begin{example}
\label{ex:psatzexample}
 Consider the polynomial system $\{f = 0, \, g \geq 0\}$, where
 \[
 f \, := \, x_2+x_1^2+2 \, = \, 0, \qquad
 g \,:=\, x_1-x_2^2+3 \, \geq \, 0.
 \]
 By the Positivstellensatz, there are no solutions $\,( x_1,x_2) \in \R^2\,$ if and
only if there exist polynomials $t_1, s_1, s_2 \in \R[x_1,x_2]$ that
 satisfy
\begin{equation}
s_1 + s_2 \cdot g +  t_1 \cdot f \,\equiv \, -1 \, , \quad
  \hbox{where $s_1$ and $s_2$ are SOS}.
 \label{eq:pstz}
 \end{equation}

 At the $D$-th SDP relaxation of the polynomial problem $\{f = 0, \, g
 \geq 0\}$, one asks whether there exists a solution $(t_1, s_1,s_2)$
 to (\ref{eq:pstz}) where the polynomial $s_1$ has degree $\leq D$ and
 the polynomials $s_2, t_1$ have degree $\leq D-2$. For each fixed
 positive integer $D$ this can be tested by a (possibly large)
 semidefinite program. Solving this for $D=2$, we find the
 infeasibility certificate
\[
 {\textstyle 
 s_1 = \frac{1}{3} 
 + 2 \left( x_2+\frac{3}{2}\right)^2
 + 6 \left( x_1-\frac{1}{6}\right)^2}, \qquad
 s_2 = 2, \qquad
 t_1 = -6.
 \]
 The resulting identity (\ref{eq:pstz}) proves the inconsistency of
 the system.
 \end{example}

As outlined in the preceding paragraphs, there is a direct connection
going from general polynomial optimization problems to SDP, via the
Positivstellensatz infeasibility certificates.  Even though we have
discussed only feasibility problems here, there are obvious
straightforward connections with optimization. For instance, by
considering the emptiness of the sublevel sets of the objective
function, or using representation theorems for positive polynomials,
sequences of converging bounds indexed by certificate degree can be
directly constructed; see e.g. \cite{Parrilo2000, Lasserre2000,
ParriloSturmfels2003}.  These schemes have been implemented in
software packages such as SOSTOOLS \cite{SOSTOOLS}, GloptiPoly
\cite{GLOPTIPOLY}, and YALMIP \cite{YALMIP}.

\subsection{Semidefinite programming relaxations}

In the last section, we have described the search for
Positivstellensatz infeasibility certificates formulated as a
semidefinite programming problem. We now describe an alternative
interpretation, obtained by dualizing the corresponding semidefinite
programs. This is the exact analogue of the construction presented in
Section~\ref{sec:LArelax}, and is closely related to the approach via
truncated moment sequences developed by Lasserre \cite{Lasserre2000}.

Recall that in the approach in Section~\ref{sec:LArelax}, the linear
relaxations were constructed by replacing every monomial $x^\alpha$ by
a new variable $\lambda_{x^\alpha}$. Furthermore, new redundant
equations were obtained by multiplying an existing constraint $f(x)=0$
by terms of the form $x_i$, yielding $x_i f(x)=0$ (essentially,
generating the ideal of valid equations). In the inequality case, and
as suggested by the Positivstellensatz, new inequality constraints
will be generated by both squarefree multiplication of the original
constraints, and by multiplication against sums of squares. That is,
if $g_i(x) \geq 0$ and $g_j(x) \geq 0$ are valid inequalities, then so
are $g_i(x) g_j(x) \geq 0$ and $g_i(x) s(x) \geq 0$, where $s(x)$ is
SOS. After substitution with the extended variables $\lambda$, we then
obtain a new system of linear equations and inequalities, with the
property that the resulting inequality conditions are
\emph{semidefinite} conditions. The presence of the semidefinite
constraints arises because we do not specify \emph{a priori} what the
multipliers $s(x)$ are, but only give their linear span.

\begin{example}
Consider the problem discussed earlier in
Example~\ref{ex:psatzexample}. The corresponding relaxation is (for
$D=2$):
\[
\begin{bmatrix}
\lambda_1 & \lambda_{x_1} & \lambda_{x_2} \\
\lambda_{x_1} & \lambda_{x_1^2} & \lambda_{x_1 x_2} \\
\lambda_{x_2} & \lambda_{x_1 x_2} & \lambda_{x_2^2} 
\end{bmatrix}\succeq 0, \qquad
\lambda_{x_2}+\lambda_{x_1^2}+ 2 \lambda_1 = 0, \qquad
\lambda_{x_1}-\lambda_{x_2^2}+3 \lambda_1 \geq 0,
\]
plus the condition $\lambda_1 > 0$ (without loss of generality, we can
take $\lambda_1 =1$). This is a semidefinite programming problem, and
in this case, its infeasibility directly shows that the original
system of polynomial inequalities does not have a solution.
\end{example}

An appealing geometric interpretation follows from considering the
projection of the feasible set of these relaxations in the space of
original variables (i.e., $\lambda_{x_i}$). For the linear algebra
relaxations of Section~\ref{sec:LArelax}, we obtain outer
approximations to the \emph{affine hull} of the solution set (an
algebraic variety), while the SDP relaxation described here constructs
outer approximations to the \emph{convex hull} of the corresponding
semialgebraic set. This latter viewpoint will be further discussed in
Section~\ref{sec:theta}, for the case of equations arising from
combinatorial problems.

\subsection{Theta bodies}
\label{sec:theta}
Recall that traditional modeling of combinatorial optimization
problems often uses $0/1$ incidence vectors. The set $S$ of solutions
of a combinatorial problem (e.g., the stable sets, traveling salesman
tours) is often computed through the (implicit) convex hull of such
incidence vectors.  Just as in the stable set and max-cut examples in
Proposition~\ref{prop:encodings}, the incidence vectors can be seen at
the set of {\em real} solutions to a system of polynomial equations:
$f_1(x) = f_2(x) = \cdots = f_m(x) = 0$, where $f_1, \ldots, f_m \in
\R[x] := \R[x_1,\ldots,x_n]$. Over the years there have been
well-known attempts to understand the structure of these convex hulls
through semidefinite programming relaxations (see
\cite{SheraliAdams1990, LovaszSchrijver1991, Lasserre2002,
Lovasz2003}) and in fact they are closely related \cite{Laurent2003,
Laurent2009}. Here we wish to summarize some recent results that give
appealing structural properties, in terms of the associated system of
equations (see \cite{GouveiaParriloThomas2008, newgouveiaetal2009} for
details).

Let us start with a historically important example: Given an
undirected finite graph $G = (V, E)$, consider the set $S_G $ of
characteristic vectors of stable sets of $G$. The convex hull of
$S_G$, denoted by $\textup{STAB}(G)$, is the \emph{stable set
polytope}.  As we mentioned already the vanishing ideal of $S_G$ is
given by $I_G := \langle x_i^2 - x_i \,\,(\forall \,\,i \in V),
\,\,\,x_ix_j \,\,(\forall \,\, \{i,j\} \in E) \rangle$ which is a real
radical zero-dimensional ideal in $\R[x]$.  In \cite{Lovasz1994},
Lov{\'a}sz introduced a semidefinite relaxation, $\textup{TH}(G)$, of
the polytope $\textup{STAB}(G)$, called the {\em theta body} of $G$.
There are multiple descriptions of $\textup{TH}(G)$, but the one in
\cite[Lemma 2.17]{LovaszSchrijver1991}, for instance, shows that
$\textup{TH}(G)$ can be defined completely in terms of the polynomial
system $I_G$. It is easy to show that $\textup{STAB}(G) \subseteq
\textup{TH}(G)$, and remarkably, we have that $\textup{STAB}(G) =
\textup{TH}(G)$ if and only if the graph is \emph{perfect}. We will
now explain how the case of stable sets can be generalized to
construct theta bodies for many other combinatorial problems.

We will construct an approximation of the convex hull of a finite set
of points $S$, denoted $\textup{conv}(S)$, by a sequence of convex
bodies recovered from ``degree truncations'' of the defining
polynomial systems.  In what follows $I$ will be a radical polynomial
ideal. A polynomial $f$ is {\bf non-negative} modulo $I$, written as
$f \geq 0 \,\,\textup{mod}\,\,I$, if $f(s) \geq 0$ for all $s \in
\Var{\R}{I}$. More strongly, the polynomial $f$ is a {\bf sum of
squares (sos)} mod $I$ if there exists $h_j \in \R[x]$ such that $f
\equiv \sum_{j=1}^t h_j^2 \,\,\textup{mod}\,\,I$ for some $t$, or
equivalently, $f-\sum_{j=1}^t h_j^2 \in I$. If, in addition, each
$h_j$ has degree at most $k$, then we say that $f$ is {\bf $k$-sos}
mod $I$.  The ideal $I$ is {\bf $k$-sos} if {\em every} polynomial
that is non-negative mod $I$ is $k$-sos mod $I$. If every polynomial
of degree at most $d$ that is non-negative mod $I$ is $k$-sos mod $I$,
we say that $I$ is {\bf $(d,k)$-sos}.  We say that a polynomial ideal
$I \subseteq \R[x]$ is $\textup{TH}_k$-exact if every linear
polynomial that is non-negative over $\Var{\R}{I}$, the real variety
of $I$, is a sum of squares of polynomials of degree at most $k$
modulo $I$.

Note that $\textup{conv}(\Var{\R}{I})$, the convex hull of
$\Var{\R}{I}$, is described by the linear polynomials $f$ such that $f
\geq 0$ mod $I$. A certificate for the non-negativity of $f$ mod $I$
is the existence of a sos-polynomial $\sum_{j=1}^{t} h_j^2$ that is
congruent to $f$ mod $I$. One can now investigate the convex hull of
$S$ through the hierarchy of nested closed convex sets defined by the
semidefinite programming relaxations of the set of $(d,k)$-sos
polynomials.

\begin{definition} \label{def:theta} 
Let $I\subseteq \R[x]$ be an ideal, and let $k$ be a positive integer.
Let $\Sigma_k \subset \R[x]$ be the set of all polynomials that are {\bf $k$-sos} mod I. 
\begin{enumerate}
\item The $k$-th {\bf theta body} of $I$ is
$$\textup{TH}_k(I) := \{ x \in \R^n \,:\,
  f(x) \geq 0 \,\,\textup{for every linear}
  \,\,f \in \Sigma_k \}.$$
\item The ideal $I$ is {\bf $\textup{TH}_k$-exact} if the $k$-th theta body
$\textup{TH}_k(I)$ coincides with the closure of $\textup{conv}(\Var{\R}{I})$.  
\item The {\bf theta-rank} of $I$ is the smallest $k$ such that
$\textup{TH}_k(I)$ coincides with the closure of
$\textup{conv}(\Var{\R}{I})$.
\end{enumerate}
\end{definition}

\begin{example} \label{ex:runningex}
  Consider the ideal $I = \langle x^2y - 1 \rangle \subset
  \R[x,y]$. Then $\textup{conv}(\Var{\R}{I}) = \{ (p_1,p_2) \in \R^2
  \,:\, p_2 > 0 \}$, and any linear polynomial that is
  non-negative over $\Var{\R}{I}$ is of the form $\alpha + \beta y$,
  where $\alpha, \beta \geq 0$. Since $\alpha y + \beta \equiv
  (\sqrt{\alpha} xy)^2 + (\sqrt{\beta})^2$ mod $I$, $I$ is $(1,2)$-sos
  and $\textup{TH}_2$-exact.
\end{example}

\begin{example}
For the case of the stable sets of a graph $G$, one can see that 
$$\textup{TH}_1(I_G) = \left\{ y \in \R^n \,:\,
\begin{array}{l} 
\exists \, M \succeq 0, M \in \R^{(n+1) \times (n+1)}\,\textup{such that} \\ 
M_{00} = 1,\\
M_{0i} = M_{i0} = M_{ii} =  y_i \,\,\forall\,\,i \in V\\
M_{ij} = 0 \,\,\forall \,\, \{i,j\} \in E
\end{array}
\right \}.$$ 
It is known that $\textup{TH}_1(I_G)$ is precisely {\em
Lov{\'a}sz's theta body of $G$}. The ideal $I_G$ is
$\textup{TH}_1$-exact precisely when the graph $G$ is perfect.
\end{example}

By definition, $\textup{TH}_1(I) \supseteq \textup{TH}_2(I) \supseteq
\cdots \supseteq \textup{conv}(\Var{\R}{I})$. As seen in
Example~\ref{ex:runningex}, $\textup{conv}(\Var{\R}{I})$ may not always
be closed and so the theta-body sequence of $I$ can converge, if at
all, only to the closure of $\textup{conv}(\Var{\R}{I})$.  But the good news
for combinatorial optimization is that there is plenty of good behavior for
problems arising with a finite set of possible solutions.

\subsection{Application: cuts and exact finite sets}
\label{sec:chull}

We discuss now a few important combinatorial examples. As we have seen
in Section~\ref{Application: Non-3-colorable graphs} for
3-colorability, and in the preceding section for stable sets, in some
special cases it is possible to give nice combinatorial
characterizations of when low-degree certificates can exactly
recognize infeasibility.  Here are a few additional results for the
real case:

\begin{example}
For the max-cut problem we saw earlier, the defining vanishing ideal
is $I(SG) = \langle x_e^2 - x_e \,\,\forall \,\,e \in E, \,\,\,
x^T \,\forall \,\, T \,\, \textup{an odd cycle in} \,\,G \rangle$.
In this case one can prove that the ideal $I(SG)$ is
$\textup{TH}_1$-exact if and only if $G$ is a bipartite graph. In
general the theta-rank of $I(SG)$ is bounded above by the size of the
max-cut in $G$. There is no constant $k$ such that
$\textup{TH}_k(I(SG)) = \textup{conv}(SG)$, for all graphs $G$. Other
formulations of max-cut are studied in \cite{newgouveiaetal2009}.
\end{example}

Recall that when $S \subset \R^n$ is a finite set, its vanishing ideal $I(S)$ is
zero-dimensional and real radical. In what follows, we say that a finite set 
  $S \subset \R^n$ is {\it exact} if its vanishing ideal $I(S)
  \subseteq \R[x]$ is $\textup{TH}_1$-exact. 

\begin{theorem}[\cite{GouveiaParriloThomas2008}] 
\label{thm:perfect}
  For a finite set $S \subset \R^n$, the following are equivalent.
\begin{enumerate} 
\item $S$ is exact.
\item There is a finite linear inequality description of
  $\textup{conv}(S)$ in which for every inequality $g(x) \geq 0$,
  $g$ is $1$-sos mod $I(S)$.
\item There is a finite linear inequality description of
  $\textup{conv}(S)$ such that for every inequality $g(x) \geq 0$,
  every point in $S$ lies either on the hyperplane $g(x)=0$ or on a
  unique parallel translate of it.
\item The polytope $\textup{conv}(S)$ is affinely equivalent to a
compressed lattice polytope (every reverse lexicographic triangulation
of the polytope is unimodular with respect to the defining lattice).
\end{enumerate}
\end{theorem}

\begin{example} \label{ex:perfect}
  The vertices of the following $0/1$-polytopes in $\R^n$ are exact
  for every $n$: (1) hypercubes, (2) (regular) cross polytopes, (3)
  hypersimplices (includes simplices), (4) joins of $2$-level
  polytopes, and (5) stable set polytopes of perfect graphs on $n$
  vertices.
\end{example}

More strongly one can say the following.

\begin{proposition} \label{prop:k parallel translates}
  Suppose $S \subseteq \R^n$ is a finite point set such that for each
  facet $F$ of $\textup{conv}(S)$ there is a hyperplane $H_F$ such
  that $H_F \cap \textup{conv}(S)=F$ and $S$ is contained in at most
  $t+1$ parallel translates of $H_F$. Then $I(S)$ is
  $\textup{TH}_t$-exact. 
\end{proposition}
 
In \cite{GouveiaParriloThomas2008} the authors show that theta bodies
can be computed explicitly as projections to the feasible set of a
semidefinite program. These SDPs are constructed using the {\em
combinatorial moment matrices} introduced by \cite{Laurent2007}.

\section{Recovering solutions in the feasible case}
\label{sec:solutions}

In principle, it is possible to find the actual roots of the system of
equations (and thus the colorings, stable sets, or desired
combinatorial object) whenever the relaxations are feasible and a few
additional conditions are satisfied. Here we discuss mostly the linear
algebra relaxations case, but the semidefinite case is very similar;
see e.g.  \cite{HenrionLasserreSolutions,LasserreLaurentRostalski2009}
for this case.

We describe below how, under certain conditions, it is possible to
recover the solution of the original polynomial system from the
relaxations (linear or semidefinite) described in earlier
sections. The main concepts are very similar for both methodologies,
and are based on the well-known eigenvalue methods for polynomial
equations; see e.g.\ \cite[\S2.4]{CoxLittleOShea05}. The key idea for
extracting solutions is the fact that from the relaxations one can
obtain a finite-dimensional representation of the vector space $R/I$
and its multiplicative structure, where $I$ is the ideal $\genR{F}$
(in the case of linear relaxations). In order to do this, we need to
compute a basis of the vector space $R/I$, and construct matrix
representations for the multiplication operators $M_{x_i}:f \mapsto
x_i f$. Then, we can use the eigenvalue/eigenvector methods to compute
solutions (see e.g., \cite{DickensteinEmiris05}).

A sufficient condition for the existence of a suitable basis of $R/I$
is given by Theorem~\ref{thm: main}. Under this condition,
multiplication matrices $M_{x_i}$ can be easily computed. In
particular, if we have computed a set $F\subset R$ that satisfies the
conditions of Theorem~\ref{thm: main} by running FPNulLA, then finding
a basis of $R/I$ and computing its multiplicative structure is
straightforward using linear algebra (see e.g., \cite{Mourrain99}). By
construction, the matrices $M_{x_i}$ commute pairwise, and to obtain
the roots one must diagonalize the corresponding commutative
algebra. It is well-known (see, e.g., \cite{CoxLittleOShea05}), that
this can be achieved by forming a random linear combination of these
matrices. This random matrix will generically have distinct
eigenvalues, and the corresponding matrix of eigenvectors will give
the needed change of basis. In the case of a finite field, it is
enough to choose the random coefficients over an algebraic extension
of sufficiently large degree, instead of working over the algebraic
closure (alternatively, the more efficient methods in
\cite{eberly2000efficient} can be used). The entries of the
diagonalized matrices directly provide the coordinates of the roots.

\begin{remark}
The condition in Theorem~(\ref{thm: main}) can in general be a strong
requirement for recovery of solutions, since it implies that we can
obtain \emph{all} solutions of the polynomial system. In some
occasions, it may be desirable to obtain just a single solution, in
which case weaker conditions may be of interest.
\end{remark}

\begin{figure}
\begin{center}
\includegraphics[height=3cm]{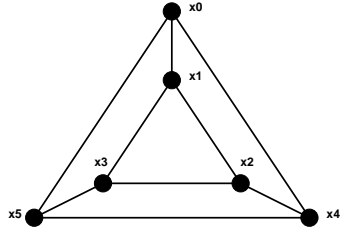}
\caption{Graph for Example~\ref{ex:extract}.}
\label{fig:twotriangles}
\end{center}
\end{figure}

\begin{example}
\label{ex:extract}
Consider the following polynomial system over $\mathbb{F}_2$, that
corresponds to the 3-colorings of the six-node graph in
Figure~\ref{fig:twotriangles}:
\[
x_i^3+1=0 \quad \forall i \in V, \qquad
x_i^2+x_i x_j + x_j^2=0 \quad \forall (i,j) \in E. 
\]
We add to these equations the symmetry-breaking constraint $x_0=1$.
After running NuLLA with this system as an input, we obtain
multiplication matrices over $\mathbb{F}_2$, of dimensions $4 \times
4$, given by:
\begin{equation*}
\begin{aligned}
M_{x_1} &= \begin{bmatrix} 0&0&1&1\\0&0&1&0\\0&1&1&0\\1&1&0&1 \end{bmatrix} &
M_{x_2} &= \begin{bmatrix} 0&0&0&1\\0&0&1&1\\1&0&1&1\\0&1&1&0 \end{bmatrix} &
M_{x_3} &= \begin{bmatrix} 0&0&1&0\\0&0&0&1\\1&1&0&1\\1&0&1&1 \end{bmatrix} \\
M_{x_4} &= \begin{bmatrix} 1&1&0&0\\1&0&0&0\\0&0&1&1\\0&0&1&0 \end{bmatrix} & 
M_{x_5} &= \begin{bmatrix} 0&1&0&0\\1&1&0&0\\0&0&0&1\\0&0&1&1 \end{bmatrix} 
\end{aligned}
\end{equation*}
Diagonalizing the corresponding commutative algebra, we obtain the
change of basis matrix given by
\[
T = \begin{bmatrix} 1&1&1&1\\
\omega^2&\omega&\omega&\omega^2\\
1&1&\omega^2&\omega\\
\omega^2&\omega&1&1 \end{bmatrix}, 
\]
where $\omega$ is a primitive root of 1, i.e., it satisfies
$w^2+w+1=0$. It can be easily verified that all the matrices $T^{-1}
M_{x_i} T$ are diagonal, and given by:
\begin{equation*}
\begin{aligned}
T^{-1} M_{x_1} T &= \diag[\omega,\omega^2,\omega,\omega^2]&
T^{-1} M_{x_2} T &= \diag[\omega^2,\omega,1,1] \\
T^{-1} M_{x_3} T &= \diag[1,1,\omega^2,\omega] &
T^{-1} M_{x_4} T &= \diag[\omega,\omega^2,\omega^2,\omega]\\
T^{-1} M_{x_5} T &= \diag[\omega^2,\omega,\omega,\omega^2],
\end{aligned}
\end{equation*}
which correspond to the four possible 3-colorings of the graph. For
instance, from the second diagonal entry of each matrix we obtain the
feasible coloring $(x_0,x_1,x_2,x_3,x_4,x_5)\rightarrow(1,\omega^2,\omega,1,\omega^2,\omega)$.
\end{example}

\bibliographystyle{siam}
\bibliography{bibliography}

\begin{thebibliography}{10}

\bibitem{BalasCeriaCornuejols1993}
{\sc E.~Balas, S.~Ceria, and G.~Cornu\'ejols}, {\em A lift-and-project cutting
  plane algorithm for mixed 0-1 programs}, Mathematical Programming, 58 (1993),
  pp.~295--324.

\bibitem{BochnakCosteRoy1998}
{\sc J.~Bochnak, M.~Coste, and M.-F. Roy}, {\em Real algebraic geometry},
  Springer, 1998.

\bibitem{CleggEdmondsImpagliazzo1996}
{\sc M.~Clegg, J.~Edmonds, and R.~Impagliazzo}, {\em Using the {G}roebner basis
  algorithm to find proofs of unsatisfiability}, in STOC '96: Proceedings of
  the twenty-eighth annual ACM symposium on Theory of computing, New York, NY,
  USA, 1996, ACM, pp.~174--183.

\bibitem{CourtoisKlimovPatarinShamir2000}
{\sc N.~Courtois, A.~Klimov, J.~Patarin, and A.~Shamir}, {\em Efficient
  algorithms for solving overdefined systems of multivariate polynomial
  equations}, in EUROCRYPT, 2000, pp.~392--407.

\bibitem{CoxLittleOShea92}
{\sc D.~Cox, J.~Little, and D.~O'Shea}, {\em Ideals, Varieties and Algorithms:
  An Introduction to Computational Algebraic Geometry and Commutative Algebra},
  Springer Verlag, 1992.

\bibitem{CoxLittleOShea05}
\leavevmode\vrule height 2pt depth -1.6pt width 23pt, {\em Using Algebraic
  Geometry}, vol.~185 of Graduate Texts in Mathematics, Springer, 2nd~ed.,
  2005.

\bibitem{DeloeraHillarMalkinWoo09}
{\sc J.~{De Loera}, C.~Hillar, P.~Malkin, and A.~Woo}, {\em Algebraic
  certificates of combinatorial infeasibility and their algorithmic
  consequences}, manuscript,  (2009).

\bibitem{DeLoeraLeeMalkinMargulies08}
{\sc J.~{De Loera}, J.~Lee, P.~Malkin, and S.~Margulies}, {\em Hilbert's
  {N}ullstellensatz and an algorithm for proving combinatorial infeasibility},
  in Proceedings of the Twenty-first International Symposium on Symbolic and
  Algebraic Computation (ISSAC 2008), 2008.

\bibitem{DeLoeraLeeMarguliesOnn08}
{\sc J.~{De Loera}, J.~Lee, S.~Margulies, and S.~Onn}, {\em Expressing
  combinatorial optimization problems by systems of polynomial equations and
  the nullstellensatz}, to appear in the Journal of Combinatorics, Probability
  and Computing,  (2008).

\bibitem{DickensteinEmiris05}
{\sc A.~Dickenstein and I.~Emiris}, eds., {\em Solving Polynomial Equations:
  Foundations, Algorithms, and Applications}, vol.~14 of Algorithms and
  Computation in Mathematics, Springer Verlag, Heidelberg, 2005.

\bibitem{eberly2000efficient}
{\sc W.~Eberly and M.~Giesbrecht}, {\em Efficient decomposition of associative
  algebras over finite fields}, Journal of Symbolic Computation, 29 (2000),
  pp.~441--458.

\bibitem{VanGelder2008}
{\sc A.~V. Gelder}, {\em Another look at graph coloring via propositional
  satisfiability}, Discrete Appl. Math., 156 (2008), pp.~230--243.

\bibitem{Gilbert1959}
{\sc E.~Gilbert}, {\em Random graphs}, Annals of Mathematical Statistics, 30
  (1959), pp.~1141--–1144.

\bibitem{newgouveiaetal2009}
{\sc J.~Gouveia, M.~Laurent, P.~A. Parrilo, and R.~R. Thomas}, {\em A new
  semidefinite programming relaxation for cycles in binary matroids and cuts in
  graphs}.
\newblock http://www.arxiv.org:0907.4518, 2009.

\bibitem{GouveiaParriloThomas2008}
{\sc J.~Gouveia, P.~A. Parrilo, and R.~R. Thomas}, {\em Theta bodies for
  polynomial ideals}.
\newblock http://www.arxiv.org:0809.3480, 2008.

\bibitem{GrigorievVorobjov2002}
{\sc D.~Grigoriev and N.~Vorobjov}, {\em Complexity of {N}ullstellensatz and
  {P}ositivstellensatz proofs}, Annals of Pure and Applied Logic, 113 (2002),
  pp.~153--160.

\bibitem{GLOPTIPOLY}
{\sc D.~Henrion and J.-B. Lasserre}, {\em {GloptiPoly}: Global optimization
  over polynomials with {MATLAB} and {SeDuMi}}, ACM Trans. Math. Softw., 29
  (2003), pp.~165--194.

\bibitem{HenrionLasserreSolutions}
\leavevmode\vrule height 2pt depth -1.6pt width 23pt, {\em Detecting global
  optimality and extracting solutions in {G}lopti{P}oly}, in Positive
  polynomials in control, vol.~312 of Lecture Notes in Control and Inform.
  Sci., Springer, Berlin, 2005, pp.~293--310.

\bibitem{HoggWilliams1994}
{\sc T.~Hogg and C.~Williams}, {\em The hardest constraint problems: a double
  phase transition}, Artif. Intell., 69 (1994), pp.~359--377.

\bibitem{KehreinKreuzer05}
{\sc A.~Kehrein and M.~Kreuzer}, {\em Characterizations of border bases},
  Journal of Pure and Applied Algebra, 196 (2005), pp.~251 -- 270.

\bibitem{KehreinKreuzerRobbiano05}
{\sc A.~Kehrein, M.~Kreuzer, and L.~Robbiano}, {\em An algebraist's view on
  border bases}, in Solving Polynomial Equations: Foundations, Algorithms, and
  Applications, A.~Dickenstein and I.~Emiris, eds., vol.~14 of Algorithms and
  Computation in Mathematics, Springer Verlag, Heidelberg, 2005, ch.~4,
  pp.~160--202.

\bibitem{Kollar1988}
{\sc J.~Koll\'ar}, {\em Sharp effective {N}ullstellensatz}, Journal of the AMS,
  1 (1988), pp.~963--–975.

\bibitem{Lasserre2000}
{\sc J.~Lasserre}, {\em Global optimization with polynomials and the problem of
  moments}, SIAM J. on Optimization, 11 (2001), pp.~796--817.

\bibitem{LasserreLaurentRostalski2009}
{\sc J.~Lasserre, M.~Laurent, and P.~Rostalski}, {\em A unified approach to
  computing real and complex zeros of zero-dimensional ideals}, in Emerging
  Applications of Algebraic Geometry, M.~Putinar and S.~Sullivant, eds.,
  vol.~149 of IMA Volumes in Mathematics and its Applications, Springer, 2009,
  pp.~125--155.

\bibitem{Lasserre2002}
{\sc J.~B. Lasserre}, {\em An explicit equivalent positive semidefinite program
  for nonlinear 0-1 programs}, SIAM J. on Optimization, 12 (2002),
  pp.~756--769.

\bibitem{Laurent2003}
{\sc M.~Laurent}, {\em A comparison of the {S}herali-{A}dams,
  {L}ov\'{a}sz-{S}chrijver, and {L}asserre relaxations for 0--1 programming},
  Math. Oper. Res., 28 (2003), pp.~470--496.

\bibitem{Laurent2004}
\leavevmode\vrule height 2pt depth -1.6pt width 23pt, {\em Semidefinite
  relaxations for max-cut}, in The Sharpest Cut: The Impact of Manfred Padberg
  and His Work, M.~Gr\"otschel, ed., vol.~4 of MPS-SIAM Series in Optimization,
  SIAM, 2004, pp.~257--290.

\bibitem{Laurent2007}
\leavevmode\vrule height 2pt depth -1.6pt width 23pt, {\em Semidefinite
  representations for finite varieties}, Mathematical Programming, 109 (2007),
  pp.~1--26.

\bibitem{Laurent2009}
\leavevmode\vrule height 2pt depth -1.6pt width 23pt, {\em Sums of squares,
  moment matrices and optimization over polynomials}, in Emerging Applications
  of Algebraic Geometry, M.~Putinar and S.~Sullivant, eds., vol.~149 of IMA
  Volumes in Mathematics and its Applications, Springer, 2009, pp.~157--270.

\bibitem{YALMIP}
{\sc J.~L{\"o}fberg}, {\em {YALMIP}: A toolbox for modeling and optimization in
  {MATLAB}}, in Proceedings of the CACSD Conference, Taipei, Taiwan, 2004.

\bibitem{Lovasz1994}
{\sc L.~Lov\'{a}sz}, {\em Stable sets and polynomials}, Discrete Math., 124
  (1994), pp.~137--153.

\bibitem{Lovasz2003}
\leavevmode\vrule height 2pt depth -1.6pt width 23pt, {\em Semidefinite
  programs and combinatorial optimization}, in Recent advances in algorithms
  and combinatorics, B.~Reed and C.~Sales, eds., vol.~11 of CMS Books in
  Mathematics, Spring, New York, 2003, pp.~137--194.

\bibitem{LovaszSchrijver1991}
{\sc L.~Lov\'{a}sz and A.~Schrijver}, {\em Cones of matrices and set-functions
  and 0-1 optimization}, SIAM J. Optim., 1 (1991), pp.~166--190.

\bibitem{SusanPhd}
{\sc S.~Margulies}, {\em Computer Algebra, Combinatorics, and Complexity:
  Hilbert's Nullstellensatz and NP-Complete Problems}, PhD thesis, UC Davis,
  2008.

\bibitem{Marshall}
{\sc M.~Marshall}, {\em Positive polynomials and sums of squares}, vol.~146 of
  Mathematical surveys and monographs, American Math. Society, 2008.

\bibitem{Mourrain99}
{\sc B.~Mourrain}, {\em A new criterion for normal form algorithms}, in Proc.
  AAECC, vol.~1719 of LNCS, Springer, 1999, pp.~430--443.

\bibitem{MourrainTrebuchet08}
{\sc B.~Mourrain and P.~Tr\'ebuchet}, {\em Stable normal forms for polynomial
  system solving}, Theoretical Computer Science, 409 (2008), pp.~229 -- 240.
\newblock Symbolic-Numerical Computations.

\bibitem{NesterovSquared}
{\sc Y.~Nesterov}, {\em Squared functional systems and optimization problems},
  in High Performance Optimization, J.~F. et~al eds., ed., Kluwer Academic,
  2000, pp.~405--440.

\bibitem{Parrilo2000}
{\sc P.~A. Parrilo}, {\em Structured semidefinite programs and semialgebraic
  geometry methods in robustness and optimization}, PhD thesis, California
  Institute of Technology, May 2000.

\bibitem{Parrilo2003}
\leavevmode\vrule height 2pt depth -1.6pt width 23pt, {\em Semidefinite
  programming relaxations for semialgebraic problems}, Mathematical
  Programming, 96 (2003), pp.~293--320.

\bibitem{ParriloSturmfels2003}
{\sc P.~A. Parrilo and B.~Sturmfels}, {\em Minimizing polynomial functions}, in
  Proceedings of the DIMACS Workshop on Algorithmic and Quantitative Aspects of
  Real Algebraic Geometry in Mathematics and Computer Science (March 2001),
  S.~Basu and L.~Gonzalez-Vega, eds., American Mathematical Society, Providence
  RI, 2003, pp.~83--100.

\bibitem{SOSTOOLS}
{\sc S.~Prajna, A.~Papachristodoulou, P.~Seiler, and P.~A. Parrilo}, {\em
  {SOSTOOLS}: Sum of squares optimization toolbox for {MATLAB}}, 2004.

\bibitem{ReidZhi08}
{\sc G.~Reid and L.~Zhi}, {\em Solving polynomial systems via symbolic-numeric
  reduction to geometric involutive form}, Journal of Symbolic Computation, In
  Press, Corrected Proof (2008).

\bibitem{SchrijverLIP}
{\sc A.~Schrijver}, {\em Theory of linear and integer programming}, Wiley,
  1986.

\bibitem{SheraliAdams1990}
{\sc H.~Sherali and W.~Adams}, {\em A hierarchy of relaxations between the
  continuous and convex hull representations for zero-one programming
  problems}, SIAM Journal on Discrete Mathematics, 3 (1990), pp.~411--430.

\bibitem{Shor}
{\sc N.~Z. Shor}, {\em Class of global minimum bounds of polynomial functions},
  Cybernetics, 23 (1987), pp.~731--734.

\bibitem{Stengle72}
{\sc G.~Stengle}, {\em A {N}ullstellensatz and a {P}ositivstellensatz in
  semialgebraic geometry}, Mathematische Annalen, 207 (1973), pp.~87--97.

\bibitem{Stetter04}
{\sc H.~Stetter}, {\em Numerical Polynomial Algebra}, SIAM, 2004.

\bibitem{VandenbergheBoyd1996}
{\sc L.~Vandenberghe and S.~Boyd}, {\em Semidefinite programming}, SIAM Review,
  38 (1996), pp.~49--95.

\bibitem{zchaff}
{\sc L.~Zhang}, {\em zchaff v2007.3.12}.
\newblock Available at http://www.princeton.edu/~chaff/zchaff.html, 2007.

\end{thebibliography}
\end{document}